\documentclass{sigcomm-alternate}
\usepackage{amsfonts,latexsym,amsmath,amstext,amssymb,verbatim,epsfig,psfrag}
\usepackage{colordvi,pstcol}
\usepackage{graphicx,psboxit}
\usepackage{psfrag}

\def\ind{{\bf 1}}
\newtheorem{remark}{Remark}
\newtheorem{theorem}{Theorem}

\begin{document}
\title{D-iteration method or how to improve Gauss-Seidel method}

\numberofauthors{1}
\author{
   \alignauthor Dohy Hong\vspace{2mm}\\
   \affaddr{Alcatel-Lucent Bell Labs}\\
   \affaddr{Route de Villejust}\\
   \affaddr{91620 Nozay, France}\\
   \email{\normalsize dohy.hong@alcatel-lucent.com}
}

\date{\today}
\maketitle

\begin{abstract}
The aim of this paper is to present the recently proposed fluid diffusion based algorithm in the general context of the matrix inversion problem associated to the Gauss-Seidel method. We explain the simple intuitions that are behind this diffusion method and how it can outperform existing methods. Then we present some theoretical problems that are associated to this representation as open research problems. We also illustrate some connected problems such as the graph transformation and the PageRank problem.
\end{abstract}
\begin{psfrags}
\section{Introduction}\label{sec:intro}
In this paper, we revisit the very well known linear algebra equation problem:
$$ A.X = B $$
where $A$ is a square matrix of size $N\times N$ and $B$ a vector of size $N$ with
unknown $X$. 
There are many known approaches to solve such an equation: Gaussian elimination,
Jacobi iteration, Gauss-Seidel
iteration, SOR (successive over-relaxation), Richardson, Krylov, Gradient method, etc
cf. \cite{Golub1996, Saad, Bagnara95aunified}).

In this paper, we propose a new iteration algorithm based on the decomposition of
matrix-vector product as a fluid diffusion model. This algorithm has been initially
proposed in the context of PageRank problem \cite{dohy}. The fluid diffusion idea
was first introduced in \cite{serge}.

\subsection{Different iterative methods}
The solution of $A.X = B$ can be solved in particular with iterative methods such as
Jacobi or Gauss-Seidel iterations, when
$A$ satisfies certain conditions (e.g. strictly diagonally dominant or symmetric and positive definite).

We recall the Jacobi iteration defined by the formula:

\begin{eqnarray*}\label{eq:jacobi}
 x_i^{(k+1)} = \frac{1}{a_{ii}}\left(b_i - \sum_{j\neq i}a_{ij}x_j^{(k)}\right), i=1, 2,..,N.
\end{eqnarray*}

and the Gauss-Seidel iteration:
\begin{eqnarray*}\label{eq:seidel}
 x_i^{(k+1)} = \frac{1}{a_{ii}}\left(b_i - \sum_{j> i}a_{ij}x_j^{(k)} - \sum_{j < i}a_{ij}x_j^{(k+1)}\right), i=1, 2,..,N.
\end{eqnarray*}

Both iterations are element-wise formula. The main difference is that in Jacobi iteration
the computation of $X^{(k+1)}$ uses only the elements of $X^{(k)}$ (similar to power iteration method
in eigenvector problem) whereas in Gauss-Seidel iteration
the computation of $x_i^{(k+1)}$ exploits the elements of $x_i^{(k+1)}$ that
have already been computed for $j<i$. And this is the main reason which explains why the Gauss-Seidel method
is generally more efficient than Jacobi iteration (when the convergence condition is satisfied).
However, unlike the Jacobi method, the Gauss-Seidel method is not adapted for a distributed computation,
because the values at each iteration are dependent on the order of the choice on $i$ and there is
not much freedom.

The approach proposed here, which we will call D-iteration (as Diffusion based iteration),
is a new improvement idea that exploits the progressive update at the vector entry level
as Gauss-Seidel with the advantage of being iteration order independent which
makes it a very interesting candidate for an asynchronous distributed computation.

\subsection{Collection vs Diffusion methods}
For the sake of simplicity and intuitive explanation,
we associate the Gauss-Seidel method to an operation of {\em collection}
(one entry of vector is updated based on the previous vector based on the incoming links).

\begin{figure}[htbp]
\centering
\includegraphics[width=6cm]{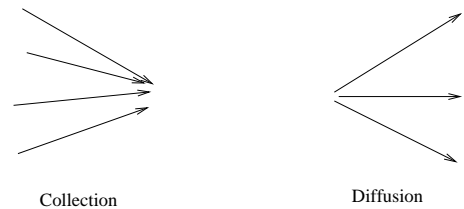}
\caption{Intuition: collection vs diffusion.}
\label{fig:cd}
\end{figure}

Our approach consists in an operation of {\em diffusion} (the fluid diffusion from one entry
of the vector consists in updating all children nodes following the outgoing links) (cf. Figure \ref{fig:cd}).
When the iteration is based on vector level update (such as Jacobi iteration or Power iteration),
the collection or diffusion approaches become equivalent (full cycle operations on all entry
of the vector).

\begin{figure}[htbp]
\centering
\includegraphics[width=\linewidth]{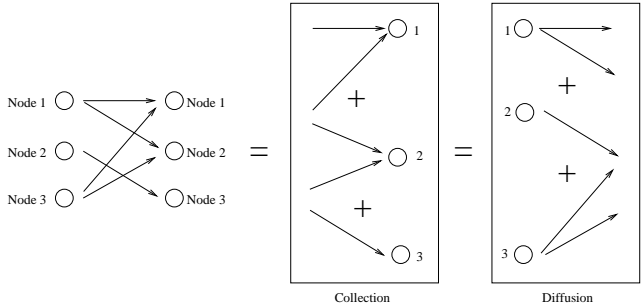}
\caption{Intuition: matrix product decomposition.}
\label{fig:eq}
\end{figure}

Somehow, those two types of operations can be seen as dual operations, but with different
consequences. With the diffusion approach, we have a very powerful result on
the convergence (monotone with an explicit information on the distance to the limit, cf. \cite{dohy})
and more importantly on the independence in the order of vector entries on which the diffusions
are applied.

\begin{figure}[htbp]
\centering
\includegraphics[width=\linewidth]{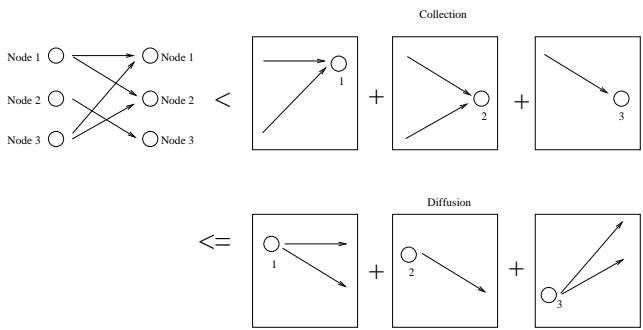}
\caption{Intuitive comparison.}
\label{fig:compa}
\end{figure}

From the intuitive point of view, D-iteration can outperform the Gauss-Seidel approach
by the appropriate optimized choice of the sequence of the vector coordinates for the diffusion
and by its distributive computation.

\section{D-algorithm}\label{sec:dalgo}
\subsection{Notations}
We reuse here the notation introduced in \cite{dohy}:
$P$ is a square matrix of size $N\times N$ such that each column sums up to 
less than one (stochastic or sub-stochastic matrix).

The D-algorithm has been initially defined in the PageRank eigenvector
context associated to the iteration of an equation of the form:
\begin{eqnarray}\label{eq:iteration0}
 X_{n+1} = A . X_n
\end{eqnarray}
where $A$ is a matrix of size $N\times N$ which can be explicitly decomposed as:
\begin{eqnarray}\label{eq:doeb}
A &=& d P + (1-d) V.\ind^t
\end{eqnarray}
where $V$ is a normalized vector of size $N$
(in the context of PageRank, the vector $V$ is a personalized initial
condition cf. \cite{deep})
and $\ind$ is the column vector with all components equal to one.
So we have:
\begin{eqnarray}\label{eq:iteration2}
 X_{n+1} = d P X_n + (1-d) V.
\end{eqnarray}

The idea of the fluid diffusion is associated to the computation of the power series:
$$
S_{\infty} = (1-d) \sum_{k=0}^{\infty} d^k P^k V = (1-d) (I_d-dP)^{-1} V = X
$$
where $I_d$ is the identity matrix and
which defines the limit of the equation \eqref{eq:iteration2}.

We define $J_k$ a matrix with all entries equal to zero except for
the $k$-th diagonal term: $(J_k)_{kk} = 1$.

In the following, we assume given a deterministic or random sequence
$I = \{i_1, i_2, ..., i_n,...\}$ with $i_n \in \{1,..,N\}$.
We only require that the number of occurrence of each value $k\in \{1,..,N\}$
in $I$ to be infinity.

We recall the definition of the two vectors used in D-iteration:
the fluid vector $F$ associated to $I$ by:
\begin{eqnarray}
F_0 &=& (1-d) V\\
F_n &=& d P J_{i_n}F_{n-1} + \sum_{k\neq i_n}J_k F_{n-1}\\
&=& (I_d - J_{i_n} + dPJ_{i_n}) F_{n-1}.\label{eq:defF}
\end{eqnarray}

And the history vector $H$ by:
\begin{eqnarray}\label{eq:defH}
H_n &=& \sum_{k=1}^n J_{i_k} F_{k-1}.
\end{eqnarray}

It is shown in \cite{dohy} that $H_n$ satisfies the iterative equation:
\begin{eqnarray*}
H_{n} &=& \left(I_d - J_{i_n}(I_d - dP)\right)H_{n-1} + J_{i_n} (1-d)V
\end{eqnarray*}
and that $H_n$ converges to the limit $X$ whatever the choice of the sequence $I$.
The $L_1$ norm of $F_n$ gives the exact distance to the limit, if the matrix $P$
has no zero column vector (otherwise, it defines an upper bound of the distance).

\subsection{Pseudo-code}
We recall here the pseudo-code presented in \cite{dohy}:

\begin{verbatim}
Initialization:
  H[i] := 0;
  F[i] := F_0[i];
  r := |F|;

Iteration:
k := 1;
While ( r/(1-d) > Target_Error )
  Choose i_k;
  sent := F[i_k];
  H[i_k] += sent;
  F[i_k] := 0;
  If ( i_k has no child )
    r -= F[i_k];
  else
    For all child node j of i_k:
      F[j] += sent*p(j,i_k)*d;
    r -= F[i_k]*(1-d);
  k++;
\end{verbatim}

$H_0$ is initialized to 0 and $F_0$ to $(1-d)V$ when associated to the equation
\ref{eq:iteration2} (the constant vector). In the case of the linear equation
$A.X = B$, $F_0$ is initialized to $B$ or $cB$ (cf. Section \ref{sec:connection}). 

\subsection{Convergence condition}
Here we only discuss the sufficient convergence condition based on the
diagonally dominant matrix.

\subsubsection{Diagonal dominant}
We recall that $A$ is strictly diagonally dominant (by columns) if:
$$
|a_{ii}| > \sum_{j\neq i} |a_{ji}|, \mbox{for all } i.
$$

For the Gauss-Seidel iteration, it is more natural to consider the row version
of the strictly diagonally dominant condition, whereas for the diffusion point
of view the column version is natural. 
However, as for the Gauss-Seidel convergence condition, both conditions guarantee
the convergence of the D-iteration.

\subsubsection{Fluid diffusion reduction}
We say that $A$ satisfies the fluid diffusion reduction condition if:
$$
\sum_{j=1}^N |a_{ji}| < 1,  \mbox{for all } i.
$$
A strictly sub-stochastic matrix (for all columns) is a specific case 
satisfying such a condition. 
This can be seen as a specific case of contractive matrix.

Here, we could also define the row version for the diffusion reduction.

Finally, as for the convergence condition of the Gauss-Seidel iteration,
the D-iteration converges if for at least one column, we have the fluid diffusion
reduction and for the other we have the equality $\sum_{j=1}^N |a_{ji}| = 1$
when $A$ is irreducible.

\subsection{Connection to the linear equation $A.X = B$}\label{sec:connection}
We can rewrite the equation $A.X = B$ as:
$$
X = (I_d - A).X + B.
$$
Then, if $X$ is a normalized probability vector $\sum_{i=1}^N x_i = 1$,
we can replace $B$ by $B.\ind^t.X$ to get:
$$
X = P.X
$$
with $P = (I_d-A) + B.\ind^t$.
We can recognize here the specific case of the PageRank equation \cite{page} for which $B= (1-d)V$
and $I_d-A = d P$ (\cite{Kohlschutter06efficientparallel, Arasu02pagerankcomputation}).

Equivalently, from an affine iteration equation
$ X_{n+1} = P X_n + B$,
we can associate to its limit $X$ a linear equation $A.X = B$ with
$A = I_d - P$.

Now we can rewrite $A.X = B$ as $cA.X = cB$ for any $c>0$ and 
$$
X = (I_d-cA).X + cB.
$$
We define $P(c) = I_d-cA$. Without any loss of generality, we can assume
that all diagonal terms of $A$ are positive (otherwise we can multiply the corresponding
line vector of $A$ and the corresponding $B$'s entry by -1).

\begin{theorem}
$A$ is strictly diagonally dominant (per column), if and only if for all $c < \frac{1}{\max_{i,j:a_{ij}\neq 0} {|a_{ij}|}}$
(and $c>0$),
$P(c)$ satisfies the fluid diffusion reduction condition.
\end{theorem}
\proof
When $A$ is strictly diagonally dominant (per column), then for all $c < \frac{1}{\max_{i,j:a_{ij}\neq 0} {|a_{ij}|}}$:
\begin{eqnarray*}
\sum_{i} |(P(c))_{ij}| &=& 1 - c|a_{ii}| + \sum_{j\neq i}c|a_{ji}| < 1.
\end{eqnarray*}
The last inequality uses exactly the strictly diagonally dominant property.
Therefore this defines a necessary and sufficient condition.

\subsection{Connection to the general eigenvector problem: $X = P.X$}\label{sec:connection2}
Here we assume that we want to solve the eigenvector equation $X = P.X$ and that we have
a square matrix $R$ such that $R.X = V$.
Then we have:
\begin{eqnarray*}
X &=& (P-R).X + V.
\end{eqnarray*}
If the spectral radius of $P-R$ is strictly less than 1, we can apply the D-iteration
to compute $X$.

If $P$ is a transition matrix and if we are looking for a probability vector $X=P.X$,
we have in particular $R = (\alpha/N) J$ where $J$ is a matrix with all entries equal to 1
and $\alpha>0$
(in case of PageRank equation, we use $R = ((1-d)/N) J$).

We illustrate this through a simple example:
take the transition matrix
$P = 
\begin{pmatrix}
0.5 & 1\\
0.5 & 0
\end{pmatrix}
$
and we want to find the stationary probability $X = P.X$.

\begin{figure}[htbp]
\centering
\includegraphics[width=\linewidth]{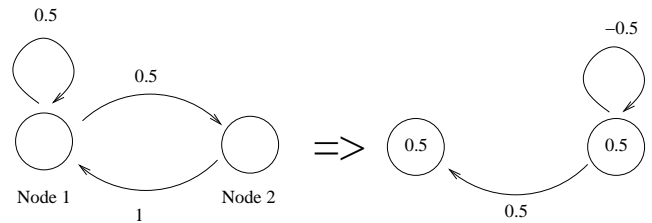}
\caption{Example: solving $X = P.X$}
\label{fig:ev}
\end{figure}

Then we have (for $\alpha = 1$):
$P-R =
\begin{pmatrix}
0 & 0.5\\
0 & -0.5
\end{pmatrix}
$

Then $X$ is can be obtained by applying the D-iteration on $(P-R,[0.5,0.5])$
or by the link elimination process defined in the following section to obtain
$X = [2/3, 1/3]$.

One simple sufficient condition on a stochastic matrix $P$ on which the D-iteration 
is convergent is given by the following theorem.

\begin{theorem}
Let $P$ be a stochastic matrix (non-negative matrix such that each column sums up to 1).
Define $N^+(i,\alpha) = |\{j: p_{ji} \ge \alpha/N\}|$ as the number of $i$-column entry of $P$ 
larger than (or equal to) $\alpha/N$.
If there exists $\alpha>0$, such that for all $i$, $N^+(i,\alpha) > N/2$, 
then the D-iteration on $(P-(\alpha/N) J, (1/N)\ind)$ is convergent.
\end{theorem}
\proof
Under the above condition, we show that $P-(\alpha/N) J$ satisfies the fluid diffusion reduction condition.
We have:
\begin{eqnarray*}
&&\sum_{i} |p_{ij} - \alpha/N| = \\
&&\sum_{i} |p_{ij} - \alpha/N|1_{p_{ij}>\alpha/N} + \sum_{i} |p_{ij} - \alpha/N|1_{p_{ij}<\alpha/N}.
\end{eqnarray*}
And we have $\sum_i p_{ij} = 1$, therefore $\sum_{i} (p_{ij} - \alpha/N) = 1 - \alpha$ and
\begin{eqnarray*}
\sum_{i} |p_{ij} - \alpha/N|1_{p_{ij}>\alpha/N} &=& 1 - \alpha + \sum_{i} |p_{ij} - \alpha/N|1_{p_{ij}<\alpha/N}.
\end{eqnarray*}
Hence
\begin{eqnarray*}
\sum_{i} |p_{ij} - \alpha/N| &=& 2 \times \sum_{i} |p_{ij} - \alpha/N|1_{p_{ij}<\alpha/N} + 1 - \alpha\\
 &\le& 2 \times \sum_{i} \alpha/N 1_{p_{ij}< \alpha/N} + 1 - \alpha\\
&\le& 2 \times (N-N^+(i,\alpha))\alpha/N + 1 - \alpha \\
&<& 1.
\end{eqnarray*}

\begin{remark}
If $P$ is irreducible, it is sufficient to have $N^+(i,\alpha) \ge N/2$
and for at least one column, a strict inequality.
\end{remark}
\begin{remark}
Obviously, the practical condition of the above theorem is to have 
$N^+(i) = |\{j: p_{ji} > 0\}| > N/2$ for all $i$.
\end{remark}

\section{Link elimination}\label{sec:trans}
The fluid diffusion model can be in the general case described
by the matrix $P$ associated with a weighted graph ($p_{ij}$ is
the weight of the edge from $j$ to $i$) and the initial condition
$F_0$. So if there is a unique limit $X$ from this setting,
we can set $X=X(P, F_0)$, which means that $X$ is the limit of the D-iteration
applied on $(P, F_0)$.

\subsection{Diagonal link elimination}
Thanks to the freedom on the sequence $I$, we have the following result:
\begin{theorem}
We have the equality (suppression of all diagonal term $p_{ii}$ such that
$p_{ii}\neq 1$; anyway if $p_{ii}\ge 1$, the D-iteration diverges):
$X(P,F_0) = X(P',F_0')$ where
\begin{itemize}
\item $(F_0')_i = (F_0)_i\times \frac{1}{1-p_{ii}}$ and $(P')_{ii} = 0$;
\item and if $i\neq j$, $(P')_{ij} = p_{ij}\frac{1}{1-p_{ii}}$.
\end{itemize}
\end{theorem}

\begin{figure}[htbp]
\centering
\includegraphics[width=\linewidth]{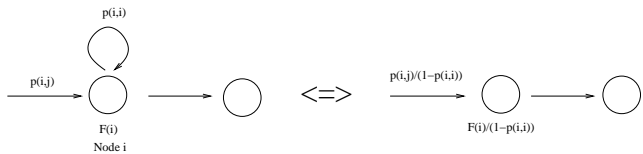}
\caption{Diagonal elimination.}
\label{fig:dig-sup-gen}
\end{figure}

\proof
The result is straightforward noticing that the self-diffusion $p_{ii}$ is
to be applied to $(F_0)_i$ and to all fluids coming from incoming links.
Such an operation is equivalent to the product by $\frac{1}{1-p_{ii}}$.




Figure \ref{fig:dig-sup} illustrates a simple example: the suppression
of the link $p_{11} = 0.75$ implies a multiplication by $1/(1-0.75) = 4$ of
$(F_0)_1$ and $p_{21}$, so that we have in this case:
$$X([[0.75, 0.5],[0.1, 0]]; [1, 1]) = X([[0, 0.5],[0.4, 0]]; [4, 1])$$.

\begin{figure}[htbp]
\centering
\includegraphics[width=\linewidth]{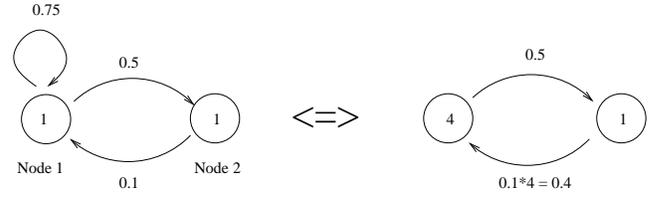}
\caption{Diagonal elimination: example.}
\label{fig:dig-sup}
\end{figure}

\subsection{Non-diagonal link elimination}
We can extend the above operation in a general case when any link
$p_{ij}$ is suppressed (but after the diagonal suppression).
\begin{theorem}
We assume that $p_{i' i'} = 0$.
Then, we have the equality (elimination of the link $p_{j'i'}$ from $i'$ to $j'$, $j'\neq i'$):
$X(P,F_0) = X(P',F_0')$ where 
\begin{itemize}
\item for $j\neq j'$, $(F_0')_j = (F_0)_j$;
\item $(F_0')_{j'} = (F_O)_{j'} + p_{j'i'}\times (F_0)_{i'}$;
\item and $(P')_{ji} = p_{ji}$ except for all $i$ an origin node of an incoming link to $i'$,
$(P')_{j'i} = p_{j'i} + p_{j'i'}\times p_{i'i}$.
\end{itemize}
\end{theorem}

\proof
The result is straightforward noticing that the elimination of the link from $i'$ to $j'$ affects
only the fluid that goes to $j'$ with the D-iteration: therefore we need to push to
$j'$ the initial fluid $(F_0)_{i'}$ and add a new link (or modify the weight of the existing one)
from all origin nodes of an incoming link to $i'$ which would replace the fluid going from
$i$ to $j'$ though $i'$.

\begin{figure}[htbp]
\centering
\includegraphics[width=\linewidth]{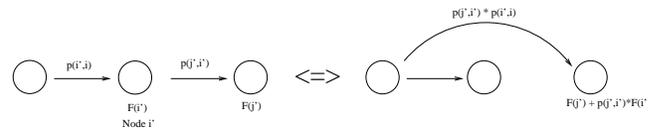}
\caption{Link $p_{j'i'}$ elimination.}
\label{fig:link-sup-gen}
\end{figure}

Figure \ref{fig:link-sup} illustrates a simple example: the suppression
of the link $p_{21} = 0.5$ implies an addition of a link $p_{22} = 0 + 0.5*0.4 = 0.2$
and the addition of $0.5\times 4 = 2$ on $(F_0)_2$ to get the equality:
$$X([[0, 0.5],[0.4, 0]]; [4, 1]) = X([[0, 0],[0.4, 0.2]]; [4, 3])$$.

\begin{figure}[htbp]
\centering
\includegraphics[width=\linewidth]{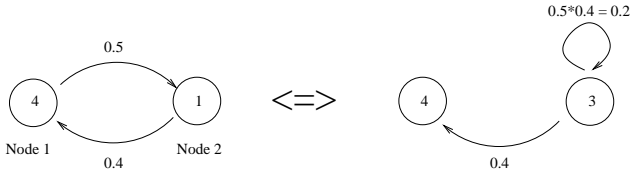}
\caption{Link $p_{21}$ elimination: example.}
\label{fig:link-sup}
\end{figure}

If we continue on suppressing the diagonal link $p_{22}$, we get:
$$X([[0, 0],[0.4, 0.2]]; [4, 3]) = X([[0, 0],[0.4, 0]]; [4, 3.75])$$.

\begin{figure}[htbp]
\centering
\includegraphics[width=\linewidth]{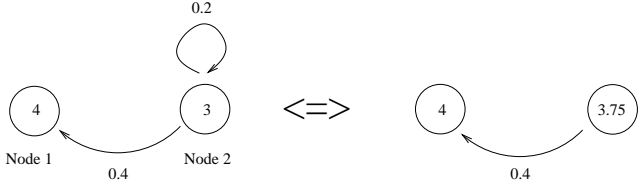}
\caption{Elimination of $p_{22}$.}
\label{fig:link-sup2}
\end{figure}

Therefore, we get:
$X = [4+0.4*3.75, 3.75] = (6, 3.75)$.

The above link suppression operation can be compared to the direct Gauss elimination method to solve
$A.X = B$ (in both cases, we get the exact limit in a finite number of operations).
The possibility of applying such a links elimination method in the computation cost reduction
of the eigenvector problem or to solve the linear equation $A.X = B$ may be an
interesting question. As the above illustration example shows, we can apply
up to the point the solution $X$ becomes explicit (no more links). However the cost of the link
elimination is proportional to the number of incoming links and in the context of
a very large matrix, such transformation may not be cost effective, because 
we may produce a very connected graph in the middle of the above transformation.

Also, it would be interesting to investigate further the idea of applying such
a transformation for the purpose of nodes clustering problem.

Below, we show one application case of such a transformation.

\subsection{Application of the graph transformation for the convergence condition}
From $P(c)$, we apply the diagonal suppression method.
The result is $Q$ such that:
$$
q_{ii} = 0
$$
and if $i \neq j$:
$$
q_{ij} = \frac{-(P(c))_{ij}}{1-(P(c))_{ii}} = \frac{-a_{ij}}{a_{ii}}.
$$
The D-iteration for $P(c)$ is convergent if and only if the D-iteration
is convergent for $Q$.
As for the Gauss-Seidel iteration, the D-iteration converges when the
spectral radius of $Q$ is strictly less than 1.

\begin{remark}
The matrix $Q$ may be obtained directly from $A.X = B$, dividing each line
$i$ by $a_{ii}$. 
\end{remark}

\subsection{Another graph transformation}
From $A.X=B$, we can apply a more natural transformation for the fluid diffusion:
because of the column based diffusion reduction condition, we set: $x_i'= x_i\times a_{ii}$
and we divide by $a_{ii}$ the $i$-th column vector of $A$ to get $A'$. Then we get
$Q' = I_d - A'$ such that:
$$
q_{ii}' = 0
$$
and if $i \neq j$:
$$
q_{ij}' = \frac{-a_{ij}}{a_{jj}}.
$$

The advantage of this approach is that this formulation is simpler to be taken into
account for the sequence $I$ optimization: for instance, for the greedy one step
vision optimization (cf. Section \ref{sec:open}).


\subsection{Example}
To illustrate the different approaches and for a simple comparison, we introduce
the following case:

$A =
\begin{pmatrix}
5 & 3 & 2 & 0\\
0 & 7 & -4 & 1\\
-2 & 0 & 8 & 0\\
0 & -2 & 1 & 3
\end{pmatrix}
$ and
$B =
\begin{pmatrix}
1\\
1\\
1\\
1
\end{pmatrix}
$ 

The results are shown on Figure \ref{fig:compa-exm}: we compared Jacobi, Power iteration
(with $c=1/8$), Gauss-Seidel, D-iteration using Q with cyclical sequence $I = {1,2,3,4,1,2 ...}$ (D-iter/Q: CYC)
and D-iteration with greedy approach taking the node with the maximum fluid in absolute value (D-iter/Q: Greedy).
For the D-iteration, we assumed here that $N$ diffusions are equivalent to one matrix-vector product iteration.
We would get the same result considering a finer iteration cost count based on the number of
link $a_{ij}$ utilization.

\begin{figure}[htbp]
\centering
\includegraphics[angle=-90,width=\linewidth]{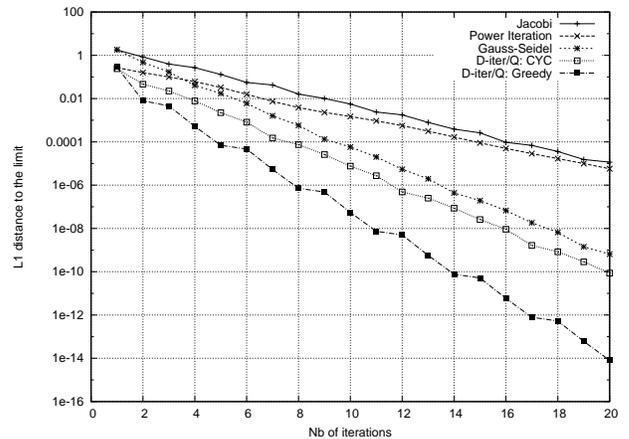}
\caption{Convergence comparison.}
\label{fig:compa-exm}
\end{figure}

To illustrate further the impact of the choice of the sequence $I$, we added
in Figure \ref{fig:compa-exm2} the result obtained when applying the D-iteration
on $Q'$ when taking the node that maximizes step by step the $L_1$-norm of $F_{n+1}$, or
equivalently the fluid reduction in one step (D-iter/Q': Greedy).

\begin{figure}[htbp]
\centering
\includegraphics[angle=-90,width=\linewidth]{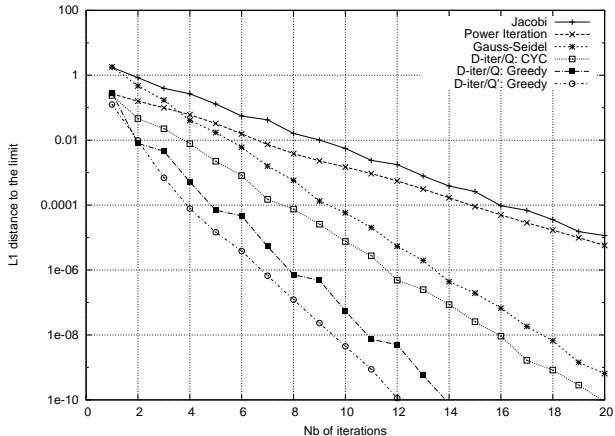}
\caption{Convergence comparison.}
\label{fig:compa-exm2}
\end{figure}

This example is only for illustration: in this particular case, the gain brought
by the D-iteration is of the order of the gain brought by the Gauss-Seidel compared
to the Jacobi iteration.

A much larger gain with D-iteration is expected with large sparse matrix
and \cite{dohy} gives an illustration of this in the context of PageRank on
the web graph. And more importantly, we can efficiently and naturally
distribute the proposed method.

\section{Open problems}\label{sec:open}


\subsection{Optimization problem}
As we showed, when the fluid reduction condition is satisfied, D-iteration
converges whatever the choice of the sequence $I$ (it is only required
that the diffusion is applied an infinity of times on each vector entry).
The convergence speed is dependent on the way the
sequence $I$ is applied. When we apply a cyclical order, the D-iteration's performance
should be intuitively close to the performance of the Gauss-Seidel iteration.
Then the natural question is: is there an optimal sequence for an optimal
convergence speed to the limit or equivalently an optimal sequence to accelerate
the fluid reduction to zero ?
An empirical solution has been proposed in \cite{dohy} in the context of PageRank associated
matrix by:
$$i_n = \arg\max_i \left(\frac{(F_{n-1})_i}{(\#in_i +1) \times (\#out_i+1)}\right).$$
In a general case when the D-iteration is applied on a matrix $Q'$, 
we could define a cost
proportional to the number of outgoing links with a gain equal to the fluid reduction factor,
so we could replace $(F_{n-1})_i/ (\#out_i+1)$ by 
$|(F_{n-1})_i|\times (1-\sum_{j}{|q_{ji}'|})/ (\#out_i+1)$.
The intuition to replace the term $1/(\#in_i +1)$ is less obvious.
So a possible empirical adaptation could be:
$$i_n = \arg\max_i \left(|(F_{n-1})_i|\times \frac{(1-\sum_{j}{|q_{ji}'|})}{(\#in_i +1)(\#out_i+1)}\right)$$
Note that we can not apply this intuitive expression to $Q$, because on one step vision, we
may have a fluid increase. Also, if $A$ is not strictly diagonally dominant, we can not use
the above expression. For the time being, in such a case, we can possibly apply again:
$$i_n = \arg\max_i \left(\frac{|(F_{n-1})_i|}{(\#in_i +1) \times (\#out_i+1)}\right).$$
The author believes that while being a very complex theoretical problem a better solution can be built
here. We hope to address those questions in a future paper. 

\subsection{Application of the graph transformation}
As it has been remarked above, it may be interesting to investigate further 
the idea of applying graph transformations described in Section \ref{sec:trans} 
for the purpose of nodes clustering problem or for a faster convergence when
mixed with the iteration methods.

\begin{remark}
It is interesting to notice that $X(P, F_0 - P.F_0) = F_0$.
This is obvious from the algebraic power series formulation of $X$, but not that
obvious from the point of view of the diffusion.
\end{remark}

\section{Conclusion}\label{sec:conclusion}
In this paper, we presented the D-iteration method, initially introduced in the
PageRank eigenvector problem, to solve efficiently
the linear equation $A.X = B$.

We believe that we have here a promising new intuitive representation and approach
that can be applied in a very large scope of linear problems.

\section*{Acknowledgments}
The author is very grateful to Fran\c cois Baccelli for taking time with patience
to follow-up this work from the beginning.
The author wishes to thank also Bruno Aidan, Paul Labrog\`ere and G\'erard Burnside
for their constant support.
Great thanks to Fabien Mathieu and St\'ephane Gaubert who put me
on the track of the Gauss-Seidel iteration.
Finally, a heartful thank to Suong Mai for her unconditional support.
\end{psfrags}
\bibliographystyle{abbrv}
\bibliography{sigproc}

\begin{thebibliography}{1}

\bibitem{serge}
S.~Abiteboul, M.~Preda, and G.~Cobena.
\newblock Adaptive on-line page importance computation.
\newblock {\em WWW2003}, pages 280--290, 2003.

\bibitem{Arasu02pagerankcomputation}
A.~Arasu, J.~Novak, J.~Tomlin, and J.~Tomlin.
\newblock Pagerank computation and the structure of the web: Experiments and
  algorithms, 2002.

\bibitem{Bagnara95aunified}
R.~Bagnara.
\newblock A unified proof for the convergence of jacobi and gauss-seidel
  methods.
\newblock {\em SIAM Review}, 37, 1995.

\bibitem{Golub1996}
G.~H. Golub and C.~F.~V. Loan.
\newblock {\em Matrix Computations}.
\newblock The Johns Hopkins University Press, 3rd edition, 1996.

\bibitem{dohy}
D.~Hong.
\newblock Optimized on-line computation of pagerank algorithm.
\newblock {\em submitted}, 2012.

\bibitem{Kohlschutter06efficientparallel}
C.~Kohlsch\"utter, P.-A. Chirita, R.~Chirita, and W.~Nejdl.
\newblock Efficient parallel computation of pagerank.
\newblock In {\em In Proc. of the 28th European Conference on Information
  Retrieval}, pages 241--252, 2006.

\bibitem{deep}
A.~N. Langville and C.~D. Meyer.
\newblock Deeper inside pagerank.
\newblock {\em Internet Mathematics}, 1(3), 2004.

\bibitem{page}
L.~Page, S.~Brin, R.~Motwani, and T.~Winograd.
\newblock The pagerank citation ranking: Bringing order to the web.
\newblock {\em Technical Report Stanford University}, 1998.

\bibitem{Saad}
Y.~Saad.
\newblock {\em Iterative Methods for Sparse Linear Systems}.
\newblock Society for Industrial and Applied Mathematics, Philadelphia, PA,
  USA, 2nd edition, 2003.

\end{thebibliography}

\end{document}